\documentclass[A4]{article}

\usepackage{graphicx}
\usepackage{url}
\usepackage{amsmath}

\usepackage{amsthm}
\usepackage{amsfonts}
\usepackage[all]{xy}
\usepackage{multicol}
\usepackage{hyperref}
\usepackage{xfrac}

\usepackage{tabularx} 
\usepackage{colortbl} 

\usepackage[percent]{overpic} 
\usepackage{subfigure} 

\usepackage{geometry}
\newcommand{\N}{\mathbb{N}}
\newcommand{\Z}{\mathbb{Z}}

\newcommand{\C}{\mathbb{C}}

\newcommand{\kbsm}{\mathcal{S}_{2,\infty}}
\newcommand{\hsm}{\mathcal{S}_3}
\newcommand{\ssm}{\mathcal{S}^{\epsilon}_{3,\infty}}
\newcommand{\sssm}{\mathcal{S}^{\pm 1}_{3,\infty}}
\newcommand{\ksm}{\mathcal{S}_{3,\infty}}
\newcommand{\dsm}{\mathcal{S}^{-1}_{3,\infty}}

\newcommand{\lpq}{L(p,q)}

\newcommand{\calS}{\mathcal{S}}\newcommand{\calL}{\mathcal{L}}

\newcommand{\rob}{\partial}

\newcommand{\lupq}{U^{-p/q}\cup L}
\newcommand{\ucl}{U \cup L}
\newcommand{\wi}{U}

\newcommand{\ov}[1]{\overline{#1}}

\newcommand{\ri}{\Omega_1} \newcommand{\rii}{\Omega_2} \newcommand{\riii}{\Omega_3}    
\newcommand{\slajd}{\textup{SL}_{p,q}}
\newcommand{\slaj}{\textup{SL}_{5,2}}
\usepackage{xcolor}

\newcommand{\mej}[1]{\text{\small{$#1$}}}

\definecolor{dark-gray}{gray}{0.3}

\newtheorem{proposition}{Proposition}

\newtheorem{theorem}{Theorem}
\newtheorem{remark}{Remark}

\newtheorem{example}{Example}

\newtheorem{conjecture}{Conjecture}

\title{Knot invariants in lens spaces}

\author{Bo\v stjan Gabrov\v sek\footnote{Faculty of Mechanical Engineering and
 Faculty of Mathematics and Physics, University of Ljubljana, Slovenia, bostjan.gabrovsek@fs.uni-lj.si}\; and
 Eva Horvat\footnote{Faculty of Education, University of Ljubljana, Slovenia, eva.horvat@pef.uni-lj.si}}

\begin{document}

\maketitle


\abstract{In this survey we summarize results regarding the Kauffman bracket, HOMFLYPT, Kauffman 2-variable and Dubrovnik skein modules, and the Alexander polynomial of links in lens spaces, which we represent as mixed link diagrams. These invariants generalize the corresponding knot polynomials in the classical case. We compare the invariants by means of the ability to distinguish between some difficult cases of knots with certain symmetries.
}

\section{Introduction} \label{sec:intro}

By the Lickorish-Wallace Theorem, any closed, connected, orientable 3-manifold $M$ can be obtained by performing Dehn surgeries on a framed link $L_0$ in $S^3$, furthermore, each component of $L_0$ can be assumed to be unknotted. Fixing $L_0$ pointwise, we can present every link $L$ in $M$ by a \emph{mixed link} $L_0 \cup L$, where we call $L_0$ the \emph{fixed component} and $L$ the \emph{moving component}, see also~\cite{lr1,dl1}.
If we take the regular projection of $L_0 \cup L$ to the plane of $L_0$, we obtain a \emph{mixed link diagram}.

\begin{figure}[ht]
	\centering\label{fig-heeg}
	\begin{overpic}[page=13]{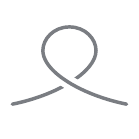}
		\put(45,-6){$m_2$}
	\end{overpic}\raisebox{1.3cm}{$\xrightarrow{\;\;h\;\;}$}
	\begin{overpic}[page=15]{images}
		\put(39,-5){$(p,-q)$}
	\end{overpic}
	\vspace{1ex}
	\caption{The boundary homeomorphism $h$.}\label{fig:heeg}
\end{figure}

In particular, if we perform $-p/q$ surgery on the unknot $U$, we obtain 
the lens space $L(p,q)$. In more detail, take $U$, remove the regular neighbourhood $\nu(U)$ of $U$ from $S^3$ and attach to the solid torus $V_1 = S^3 \setminus \nu(U)$ the solid torus $V_2 = S^1 \times D^2$ by the boundary homeomorphism $h: \partial V_2 \rightarrow \partial V_1$ that maps the meridian $m_2$ of $\partial V_2 \approx S^1 \times S^1$ to the $(p,-q)$-curve on $\partial V_1 \approx S^1 \times S^1$, which is the curve that wraps $p$-times around the longitude and $-q$-times around the meridian of $\partial V_1$ as illustrated in Figure~\ref{fig-heeg}.


\begin{figure}[ht]
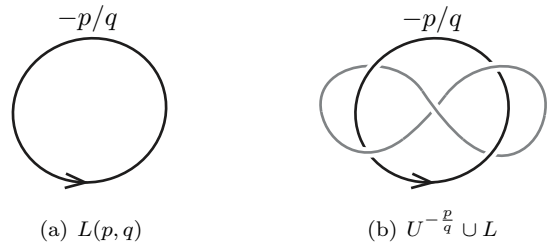

	\centering
	\subfigure[$L(p,q)$]{
	\begin{overpic}[page=46]{images}
		\put(34,85){$-p/q$}
	\end{overpic}\label{fig:kirby1}
	}\hspace{7ex}
	\subfigure[$U^{-\frac{p}{q}}\cup L$]{
	\begin{overpic}[page=45]{images}
		\put(38,70){$-p/q$}
	\end{overpic}\label{fig:kirby2}
	}
	
	\caption{The diagram of $L(p,q)$ and mixed link representing a knot in $L(p,q)$.}\label{fig:kirby}
\end{figure}

A link $L$ in $L(p,q)$ can thus be represented by the mixed link diagram of $U \cup L$. When appropriate, we will emphasize that surgery has been performed on $U$ by equipping the diagram with surgery coefficients as in Figure~\ref{fig:kirby} and we will denote such a link in $L(p,q)$ by $U^{-\frac{p}{q}}\cup L$. Note that even when dealing with unoriented links, the fixed component should be oriented, since the ambient manifold depends on this orientation.


\begin{figure}[ht]
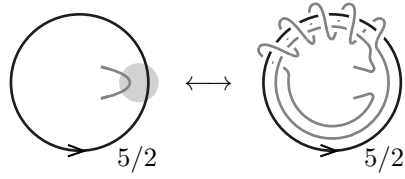

\centering

\begin{overpic}[page=66]{images}
\put(67,0){$5/2$}
\end{overpic}\raisebox{1cm}{$\;\;\;\longleftrightarrow\;$}
\begin{overpic}[page=67]{images}
\put(67,0){$5/2$}
\end{overpic}
\vspace{1ex}
\caption{The slide move $\slaj$ in $L(5,2)$.}\label{fig:slide}
\end{figure}

\begin{figure}[ht]
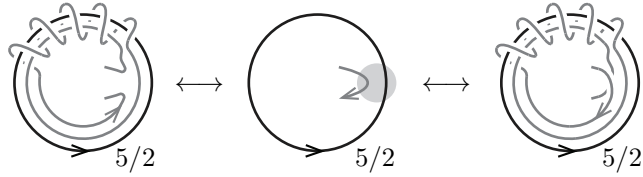

\centering
\begin{overpic}[page=64]{images}
\put(67,0){$5/2$}
\end{overpic}\raisebox{1cm}{$\;\;\longleftrightarrow\;$}
\begin{overpic}[page=62]{images}
\put(67,0){$5/2$}
\end{overpic}\raisebox{1cm}{$\;\;\;\longleftrightarrow\;$}
\begin{overpic}[page=61]{images}
\put(67,0){$5/2$}
\end{overpic}
\vspace{1ex}
\caption{Two oriented slide moves in $L(5,2)$.}\label{fig:slide2}
\end{figure}

If we approach the meridian disk of $V_2$ with an arc of $L$, we can slide the arc along the disk bounding $m_2$ (the 2-handle in the CW decomposition of $L(p,q)$), which has the effect of making a connected sum with $(p,-q)$-curve representing $\rob m_2$ on $\rob V_1$~\cite{HP1, lr2, dl1}. This isotopy move, called the \emph{slide move} (or in some literature the \emph{band move}), is illustrated in Figure~\ref{fig:slide} and we denote it by $\slajd$. If we consider oriented links, we often differentiate between two variants of the slide move, one where the curve travels along the orientation of $U$ and the other one where we travel in the opposite direction, depending on how the approaching arc is oriented with respect to the orientation of $U$. The two oriented flavours of $\slajd$ are illustrated in Figure~\ref{fig:slide2}.

\begin{figure}[htb]
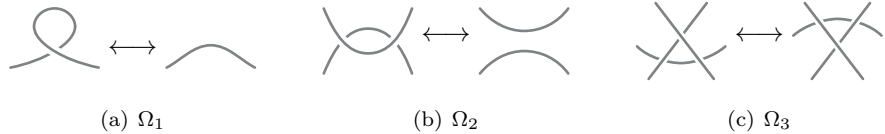

\centering
\subfigure[$\ri$]{\begin{overpic}[page=1]{images}\end{overpic}\raisebox{15pt}{$\longleftrightarrow$}\begin{overpic}[page=2]{images}\end{overpic}\label{fig:reid1}}
\hspace{1.2em}
\subfigure[$\rii$]{\begin{overpic}[page=3]{images}\end{overpic}\raisebox{20pt}{$\longleftrightarrow$}\begin{overpic}[page=4]{images}\end{overpic}}
\hspace{1.2em}
\subfigure[$\riii$]{\begin{overpic}[page=5]{images}\end{overpic}\raisebox{20pt}{$\longleftrightarrow$}\begin{overpic}[page=6]{images}\end{overpic}}
\caption{Classical Reidemeister moves.}
\label{fig:reid}
\end{figure}

The slide move, together with the planar Reidemeister moves in Figure~\ref{fig:reid} are sufficient to describe isotopy in $L(p,q)$ as the following theorem states.

\begin{theorem}[\cite{HP1}]Two mixed link diagrams represent the same link in $\lpq$ if and only if one can be transformed into the other by a finite sequence of Reidemeister moves $\ri$, $\rii$, $\riii$, and $\slajd$.\end{theorem}

\begin{remark}
Since $U$ is fixed, the arcs involved in $\ri$ belong to the moving component, in $\rii$ at most one of the arcs can belong to the fixed component and in $\riii$ at most two arcs can belong to the fixed component.\end{remark}

\section{The Kauffman Bracket skein module} \label{sec:kbsm}

Let $\raisebox{-4.5pt}{\includegraphics[page=16]{skein}},\raisebox{-4.5pt}{\includegraphics[page=17]{skein}},\raisebox{-4.5pt}{\includegraphics[page=18]{skein}}$ be the (oriented) skein triple and $\raisebox{-4.5pt}{\includegraphics[page=1]{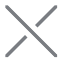}},\raisebox{-4.5pt}{\includegraphics[page=2]{skein}},\raisebox{-4.5pt}{\includegraphics[page=3]{skein}}$ the (unoriented) Kauffman triple, i.e., links that are the same everywhere except inside a small $3$-ball where they differ as the notation suggests.

Skein modules have their origin in the observation made by J. W. Alexander that the Alexander polynomials 
$\Delta\big( \raisebox{-4.5pt}{\includegraphics[page=16]{skein}}\big)$, $\Delta\big( \raisebox{-4.5pt}{\includegraphics[page=17]{skein}}\big)$, and $\Delta\big( \raisebox{-4.5pt}{\includegraphics[page=18]{skein}}\big)$ are linearly related by the skein relation
$$\Delta\big( \raisebox{-4.5pt}{\includegraphics[page=16]{skein}}\big) - \Delta\big( \raisebox{-4.5pt}{\includegraphics[page=17]{skein}}\big) = (t^{1/2} - t^{-1/2}) \Delta\big( \raisebox{-4.5pt}{\includegraphics[page=18]{skein}}\big).$$ 
J. H. Conway pursued this idea by taking $z=t^{1/2} - t^{-1/2}$ and considering the free $\Z[z]$-module over the set of isotopy classes of links in $S^3$ modulo the $\Z[z]$-module generated by the skein relation of the Alexander-Conway polynomial~\cite{kauff, Pr0, P}.

By formalizing such a construction and generalizing it for arbitrary 3-manifolds, J. H. Przytycki and V. G Turaev  introduced the theory of skein modules in~\cite{Tu, Pr2}. 

The Kauffman bracket skein module generalizes the Kauffman bracket in the following sense.

Take a coefficient ring $R$ with $A \in R$ being a unit (an element with a multiplicative inverse).
Since, as in the case of the Kauffman bracket, we would like to study framed links, we set $\calL_{\mathrm{fr}}(M)$ to be the set of isotopy classes of framed links in $M$, including the empty link $\emptyset$. Let $R\calL_{\mathrm{fr}}(M)$ be the free $R$-module spanned by $\calL_{\mathrm{fr}}(M)$.

We would like to impose the Kauffman relation and the framing relation in $R\calL_{\mathrm{fr}}(M)$. We therefore take the submodule $\mathcal{S}(M)$ of $R\calL_{\mathrm{fr}}(M)$ generated by
\begin{align*}
\tag{Kauffman relator} \raisebox{-4.5pt}{\includegraphics[page=1]{skein}} - A \raisebox{-4.5pt}{\includegraphics[page=2]{skein}} -A^{-1} \raisebox{-4.5pt}{\includegraphics[page=3]{skein}},\\
\tag{framing relator} L \sqcup \raisebox{-4.5pt}{\includegraphics[page=4]{skein}} - (-A^2 - A^{-2}) L.
\end{align*}

The \emph{Kauffman bracket skein module} $\kbsm(M)$ is $R\calL_{\mathrm{fr}}(M)$ modulo these two relations: $$\kbsm(M) = R\calL_{\mathrm{fr}}(M) / \mathcal{S}(M).$$ 

\begin{figure}[ht]
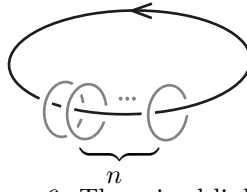

	\centering
	\begin{overpic}[page=52]{images}
		\put(30,-6){\rotatebox{90}{$\Bigg\{ $}}
		\put(40,-13){$n$}
	\end{overpic}
	\caption{The mixed link $x^n$.}\label{xn}
\end{figure}

Let $U$ be a fixed unknot in $S^3$ and let $x^n$ be the mixed link where the moving components consists of $n$ parallel copies of the unknot linked with $U$ as in Figure~\ref{xn}. Separately, we denote by $x^0$ the \emph{affine} unknot (the unknot contained inside a 3-ball in $M$).

If we remove a tubular neighbourhood $\nu(U)$ of $U$, we can think of $U\cup L$ as a link in the solid torus $T = V_1$. 

The Kauffman bracket skein module of the solid torus $T$ has been calculated by Turaev:
\begin{theorem}[Turaev~\cite{Tu}]\label{theorem:kbsmt}
$\kbsm{T}$ is a free $R$-module generated by the set $\{x^n\}_{n=0}^{\infty}$.
\end{theorem}

If, instead of removing $U$, we perform $-p/q$ surgery on $U$, we can think of $x^n$ as a link in $L(p,q)$.

\begin{theorem}[Hoste, Przytycki~\cite{HP1}]\label{theorem:kbsmlpq}
$\kbsm(\lpq)$ is a free $R$-module generated by $\{x^n\}_{n=0}^{\lfloor p/2 \rfloor}$.
\end{theorem}

These generating sets are just natural choices, for alternative bases see~\cite{thales}.
The KSBM has been calculated for several other classes of manifolds, see for example~\cite{mr0, mr1, mr2}.

\section{The HOMFLYPT skein module} \label{sec:hom}

The HOMFLYPT skein module  of a 3-manifold $M$ generalizes the HOMFLYPT polynomial. Let the ring $R$ this time have two units $v, z \in R$. Let $\calL_{\mathrm{or}}(M)$ be the set of isotopy classes of oriented links in $M$, including the empty link $\emptyset$ and let $R\calL_{\mathrm{or}}(M)$ be the free $R$-module spanned by $\calL_{\mathrm{or}}(M)$. 

We impose the HOMFLYPT skein relation in $R\calL_{\mathrm{or}}(M)$ by taking the submodule $\calS(M)$ of $R\calL_{\mathrm{or}}(M)$ generated by the expressions 
\begin{align*}
\tag{HOMFLYPT relator} v^{-1} \raisebox{-4.5pt}{\includegraphics[page=16]{skein}} - v \raisebox{-4.5pt}{\includegraphics[page=17]{skein}} - z \raisebox{-4.5pt}{\includegraphics[page=18]{skein}}.
\end{align*}

We also add to $\calS(M)$ the HOMFLYPT relation involving the empty knot,
\begin{align*}
\tag{HOMFLYPT relator} v^{-1}\emptyset - v\emptyset - z\raisebox{-4.5pt}{\includegraphics[page=20]{skein}}.
\end{align*}

The \emph{HOMFLYPT skein module} $\hsm(M)$ of $M$  is $R\calL_{\mathrm{or}}(M)$ modulo the above relations:
$$\hsm(M) = R\calL_{\mathrm{or}}(M)/\calS(M).$$

Let $U$ be a fixed unknot and let $t_k$, $k \in \Z \setminus \{0\}$, be the oriented link that wraps $k$ times around $U$ as in Figures~\ref{fig:hsm1} and~\ref{fig:hsm2}(note that $t_{-k}$ is $t_k$ with reversed orientation). We define the product $t_{k_1} t_{k_2} \cdots t_{k_s}$, $s\in \N$, as the links $t_{k_i}$ placed consecutively along $U$ as illustrated in Figure~\ref{fig:hsm3}.

\begin{figure}[ht]
	\centering
	\subfigure[$t_k, k>0$]{
	\begin{overpic}[page=55]{images}
				\put(32,43){\rotatebox{-90}{$\Bigg\{ $}}
		\put(44,44){$k$}
	\end{overpic}\label{fig:hsm1}
	}\hspace{1em}
	\subfigure[$t_k, k<0$]{
	\begin{overpic}[page=56]{images}
					\put(32,43){\rotatebox{-90}{$\Bigg\{ $}}
		\put(44,44){$k$}
	\end{overpic}\label{fig:hsm2}
	}
	\subfigure[$t_{-2}^2 t_3$]{
	\begin{overpic}[page=53]{images}
	\end{overpic}\label{fig:hsm3}
	}
	\caption{Generators of $\hsm(T)$.}\label{fig:hsm}
\end{figure}

\begin{theorem}[Turaev~\cite{Tu}]\label{theorem:hsmt}
$\hsm(T)$ is a free $R$-module generated by 
$$\{t^{i_1}_{k_1} \ldots t^{i_s}_{k_s} \mid s\in\N,\; k_j \in \Z \setminus \{0 \},\;
k_1< \cdots <k_s,\; i_j \in \N\}\cup\{\emptyset\}.$$
\end{theorem}

\begin{theorem}[\cite{gm}]\label{theorem:hsmt1}
$\hsm(L(p,1))$ is a free $R$-module generated by 
$$\{ t_{k_1}^{i_1}\cdots t_{k_s}^{i_s} \mid s\in\mathbb N,\; k_j \in \Z \setminus \{0\}, -\frac{p}{2}<k_1 < \cdots < k_s \leq \frac{p}{2},\; i_j \in \N \}.$$
\end{theorem}
For alternative bases see~\cite{thales} and~\cite{dl2}. The proof of Theorem~\ref{theorem:hsmt1} in~\cite{gm} is based on a diagramatic approach, but the problem can be also attacked using a braid approach, see~\cite{dlp, dl3}.

The case of $\hsm(L(p,q))$, $q \geq 2$, is still an open question, but it is believed that the following conjecture holds.

\begin{conjecture}\label{theorem:hsmt2}
$\hsm(L(p,q))$ is a free $R$-module generated by 
$$\{ t_{k_1}^{i_1}\cdots t_{k_s}^{i_s} \mid s\in\mathbb N,\; k_j \in \Z \setminus \{0\}, -\frac{p}{2}<k_1 < \cdots < k_s \leq \frac{p}{2},\; i_j \in \N \}.$$
\end{conjecture}

Related to this invariant, in~\cite{corn} Cornwell constructed a 2-variable polynomial in $L(p,q)$ that satisfies the skein relation (but is in essence weaker than the HOMFLYPT skein module), see also~\cite{cmr} where this invariant has been studied.

\section{The Kauffman and Dubrovnik skein modules} \label{sec:dubrovnik}

The Kauffman and Dubrovnik skein modules generalize the Kauffman 2-variable and Dubrovnik polynomials of unoriented links.

Let the ring $R$ have two units $z, a \in R$.
Take the submodule $\calS(M)$ of $R\calL_{\mathrm{fr}}(M)$ generated by the expressions 
\begin{align*}
\tag{Kauffman/Dubrovnik relator} \raisebox{-4.5pt}{\includegraphics[page=21]{skein}} +
\epsilon \raisebox{-4.5pt}{\includegraphics[page=1]{skein}}
-z \raisebox{-4.5pt}{\includegraphics[page=3]{skein}}
-\epsilon z \raisebox{-4.5pt}{\includegraphics[page=2]{skein}},\\
\tag{framing relator} \raisebox{-4.5pt}{\includegraphics[page=23]{skein}} 
- a \raisebox{-4.5pt}{\includegraphics[page=22]{skein}} .
\end{align*}

We add to $\calS(M)$ the relation involving the empty knot,
\begin{align*}
\raisebox{-4.5pt}{\includegraphics[page=4]{skein}} - \Big(\frac{u+\epsilon u^{-1}}{z}- \epsilon\Big) \emptyset.
\end{align*}

We define the module
$$\ssm(M) = R\calL(M)/\calS(M).$$

Taking $\epsilon = +1$, we obtain the \emph{Kauffman skein module} $\ksm(M)$ and for $\epsilon = -1$, we obtain the \emph{Dubrovnik skein module} $\dsm(M)$.

Let $t_k$, $k \in \N \setminus \{0\}$, be the unoriented knot that wraps $k$ times around $U$ as in Figure~\ref{fig:hsm1}. As in the previous section, the product $t_{k_1} t_{k_2} \cdots t_{k_s}$, $s\in \N$ is the link consisting 
of $t_{k_i}$'s placed along $U$ as illustrated in Figure~\ref{fig:hsm3}.


\begin{figure}[ht]
	\centering
	\subfigure[$t_k$]{
	\begin{overpic}[page=58]{images}
				\put(32,43){\rotatebox{-90}{$\Bigg\{ $}}
		\put(44,44){$k$}
	\end{overpic}\label{fig:ksm1}
	}\hspace{2em}
	\subfigure[$t_{1}^2 t_3$]{
	\begin{overpic}[page=59]{images}
	\end{overpic}\label{fig:ksm2}
	}
	\caption{Generators of $\sssm(T)$.}\label{fig:ssm}
\end{figure}

For the solid torus both modules have been calculated in by Turaev: 

\begin{theorem}[Turaev~\cite{Tu}]\label{theorem:ssmt}
$\sssm(T)$  are free $R$-modules generated by 
$$\{ t_{k_1}^{i_1}\cdots t_{k_s}^{i_s} \mid s, k_j \in \mathbb N,\; 0<k_1 < \cdots < k_s,\; i_j \in \N \} \cup \{ \emptyset \}.$$
\end{theorem}

For the lens spaces $L(p,1)$ the modules have been calculated by Mroczkowski:

\begin{theorem}[Mroczkowski~\cite{mr3}]\label{theorem:lpqsm1}
$\ksm(L(p,1))$ is generated by 
$$\{ t_{k_1}^{i_1}\cdots t_{k_s}^{i_s} \mid s, k_j \in \mathbb N,\; 0<k_1 < \cdots < k_s \leq \lfloor \frac{p}{2}\rfloor,\; i_j \in \N \} \cup \{ \emptyset \}.$$
The modules are free if $p$ is odd and contain torsion if $p$ is even. 
\end{theorem}

\begin{theorem}[Mroczkowski~\cite{mr3}]\label{theorem:lpqsm2}
$\dsm(L(p,1))$ is a free $R$-module generated by 
$$\{ t_{k_1}^{i_1}\cdots t_{k_s}^{i_s} \mid s, k_j \in \mathbb N,\; 0<k_1 < \cdots < k_s \leq \lfloor \frac{p}{2}\rfloor,\; i_j \in \N \} \cup \{ \emptyset \}.$$
\end{theorem}

\section{Alexander polynomial} \label{sec:alexander}

In this section we describe a Torres-type formula (see \cite{torres}), constructed in~\cite{gh} for the Alexander polynomial of links in lens spaces defined by Fox's free differential calculus~\cite{fox, lin, wada}.

\begin{figure}[ht]
	\centering
	\subfigure[positive crossing]{\begin{overpic}[page=1]{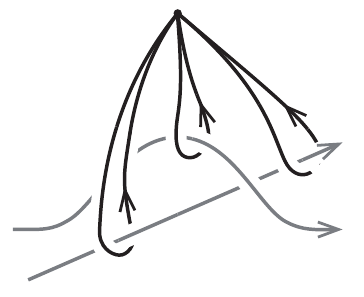}
		\put(19,44){$x_{i_1}$}\put(72,61){$x_{i_2}$}\put(44,33){$x_{i_3}$}
		\label{fig:wirt1}
	\end{overpic}}
	\qquad\qquad\qquad
	\subfigure[negative crossing]{\begin{overpic}[page=2]{slike}
		\put(19,44){$x_{i_1}$}\put(72,61){$x_{i_2}$}\put(44,33){$x_{i_3}$}
		\label{fig:wirt2}
	\end{overpic}}
	\caption{Wirtinger relations.}\label{fig:wirt}
\end{figure}
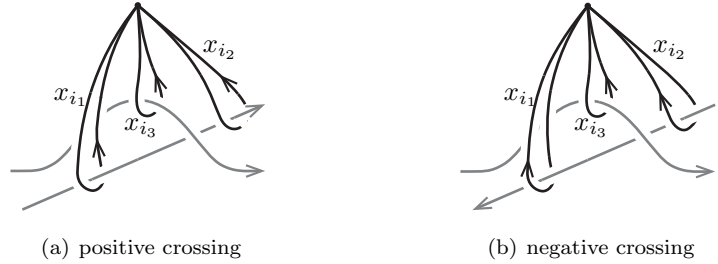

Recall that the fundamental group of a classical link admits a well-known Wirtinger presentation 
$$\pi _{1}(S^{3}\backslash L,*)=\left \langle x_{1},\ldots ,x_{n} \mid r_{1},\ldots ,r_{n}\right \rangle ,$$ 
obtained from a link diagram. Generators $x_{i}$ correspond to the simple closed loops based at $*$ and winding around the over-arcs of the diagram and $r_{i}$ is the Wirtinger relation, $x_{i_{1}}x_{i_{3}}x_{i_{2}}^{-1}x_{i_{3}}^{-1}$ if the crossing is positive or $x_{i_{1}}x_{i_{3}}^{-1}x_{i_{2}}^{-1}x_{i_{3}}$ if the crossings is negative, corresponding to the $i$-th crossing of the diagram, see Figure \ref{fig:wirt}. 

Given a mixed link diagram of $\lupq$ the following proposition allows us to describe the fundamental group of $L(p,q) \setminus L$ (cf.~\cite{cmm,gma}).

\begin{proposition}[\cite{rolfsen}]\label{prop1}
Let $\left \langle x_{1},\ldots ,x_{n} \mid \, r_{1},\ldots ,r_{n}\right \rangle $ be the Wirtinger presentation for $\pi _{1}(S^{3}\backslash (\ucl),*)$ obtained from a mixed link diagram. Denote by $m_{1}$ and $l_{1}$ the meridian and longitude of the regular neighbourhood of $S^{3}\backslash \wi$, written in terms of the generators $x_{1},\ldots ,x_{n}$. The presentation for the link group is given by 
$$ \pi _{1}(L(p,q)\backslash L,*)=\left \langle x_{1},\ldots ,x_{n} \mid \, w_{1},\ldots ,w_{n},m_{1}^{p}\, l_{1}^{-q}\right \rangle.$$
\end{proposition} 

We briefly recall the construction of the Alexander polynomial using Fox calculus \cite{wada, gh}. Suppose $$\mathcal{P}=\left \langle x_{1},\ldots ,x_{n} \mid \, r_{1},\ldots ,r_{m}\right \rangle $$ is a presentation of a group $G$. Denote by $H=G/G'$ its abelianization and by $F=\langle x_{1},\ldots ,x_{n} \mid \, \rangle $ the corresponding free group. Apply the chain of maps $$\Z F\stackrel{\frac{\partial }{\partial x}}{\longrightarrow }\Z F\stackrel{\gamma }{\longrightarrow}\Z G\stackrel{\alpha }{\longrightarrow}\Z H\;,$$ where $\frac{\partial }{\partial x}$ denotes the Fox differential, $\gamma $ is the quotient map by the relations $r_{1},\ldots ,r_{m}$ and $\alpha $ is the abelianization map. 

The \emph{Alexander-Fox matrix} of $\mathcal{P}$ is the matrix $A=\left [a_{i,j}\right ]$, where $a_{i,j}=\alpha (\gamma (\frac{\partial r_{i}}{\partial x_{j}}))$ for $i=1,\ldots ,m$ and $j=1,\ldots n$. The \emph{first elementary ideal} $E_{1}(\mathcal{P})$ is the ideal of $\Z H$, generated by the determinants of all the $(n-1)$ minors of $A$. 

For a link $L$ in $S^3$, let $E_{1}(\mathcal{P})$ be the first elementary ideal obtained from a presentation $\mathcal{P}$ of $\pi _{1}(S^{3}\backslash L,*)$. The \emph{Alexander polynomial} $\Delta (L)$ is the generator of the smallest principal ideal containing $E_{1}(\mathcal{P})$. The abelianization of $\pi _{1}(S^{3}\backslash L,*)$ is a free abelian group whose generators correspond to the components of $L$. 

For a link in $L(p,q)$, the abelianization of its link group may also contain torsion, see \cite[Corollary 2.10]{gh}. In this case, we need the notion of a twisted Alexander polynomial. We recall the following from~\cite{cmm}. 

Let $G$ be a group with a finite presentation $\mathcal{P}$ and abelianization $H=G/G'$ and denote $K=H/Tors(H)$. Then every homomorphism $\sigma \colon Tors(H)\to \C ^{*}=\C \backslash \{0\}$ determines a twisted Alexander polynomial $\Delta ^{\sigma }(\mathcal{P})$ as follows. Choosing a splitting $H=Tors(H)\times K$, $\sigma $ defines a ring homomorphism $\sigma \colon \Z [H]\to \C [G]$ sending $(f,g)\in Tors(H)\times K$ to $\sigma (f)g$. Thus we apply the chain of maps $$\Z F\stackrel{\frac{\partial }{\partial x}}{\longrightarrow }\Z F\stackrel{\gamma }{\longrightarrow}\Z G\stackrel{\alpha }{\longrightarrow}\Z H\stackrel{\sigma }{\longrightarrow}\C [K]$$ and obtain the $\sigma $-twisted Alexander matrix $A^{\sigma }=\left [\sigma (\alpha (\gamma (\frac{\partial r_{i}}{\partial x_{j}})))\right ]$. The \emph{twisted Alexander polynomial} is then defined by $\Delta ^{\sigma }(\mathcal{P})=\textrm{gcd}(\sigma (E_{1}(\mathcal{P})))$. 

The \emph{Alexander polynomial} of $\lupq$, which we denote by $\Delta_{\lupq}$ or simply $\Delta_L$ if the context is clear, is defined to be the generator of the smallest principal ideal containing $E_{1}(\mathcal{P})$. 

We continue by describing how to obtain the Alexander polynomial for $\lupq$ from the Alexander polynomial of $\ucl \subset S^3$.

Let $D$ be the disk bounded by $\wi$. We may assume that $L$ intersects $D$ transversely in $k$ intersection points with algebraic intersection signs $\epsilon _{1},\ldots ,\epsilon _{k}\in \{-1,1\}$. We define $[L] = \sum_{i=1}^k \epsilon_i$, which corresponds to the integer representing the homology class of $L$ in $H_1(S^3 \setminus \wi) \cong \Z$.

By Proposition~\ref{prop1}, the presentation of $\pi_1(L(p,q) \setminus L,*)$ is obtained from the presentation of the link group $\pi_{1}(S^{3}\setminus(\ucl),*)$ by adding one relation. The Alexander-Fox matrices are thus closely related and consequently so are the Alexander polynomials, as the following theorem states.


\begin{theorem}[\cite{gh}]
\label{th1} Let $p' = \frac{p}{\gcd\{p,[L]\}}$ and $[L]' = \begin{cases}
1, & \text{if } [L]=0\\
\frac{[L]}{\gcd\{p,[L]\}}, & \text{if } [L]\neq 0
\end{cases}.$ The Alexander polynomial of $\lupq$ and the (classical) two-variable Alexander polynomial $\Delta_{\ucl}(u,t)$, where variable $u$ corresponds to the moving components and variable $t$ corresponds to the fixed component, are related by 
\begin{equation}\label{eq:alexlpq}
\Delta_{\lupq}(t)=\frac{\Delta_{\ucl}(t^{p'},t^{q[L]'})}{t^{[L]'}-1}\;.
\end{equation}
\end{theorem}

It is also shown in~\cite{gh} that it is possible to normalize $\Delta_{\lupq}$ and obtain a normalized version of the Alexander polynomial in lens spaces, $\nabla(L)(t)$, which satisfies the skein relation 
$$\nabla \big(\raisebox{-4.5pt}{\includegraphics[page=16]{skein}}\big) -\nabla\big(\raisebox{-4.5pt}{\includegraphics[page=17]{skein}}\big)=(t^{\frac{p'}{2}}-t^{-\frac{p'}{2}})\nabla \big(\raisebox{-4.5pt}{\includegraphics[page=18]{skein}}\big).$$

This result may be compared to the skein relation for links in the projective space $L(2,1)$ obtained in \cite{huynh}: 
\begin{theorem}[Huynh, Le~\cite{huynh}] Let $\raisebox{-4.5pt}{\includegraphics[page=16]{skein}},\raisebox{-4.5pt}{\includegraphics[page=17]{skein}},\raisebox{-4.5pt}{\includegraphics[page=18]{skein}}$ be a skein triple in the projective space. If $\raisebox{-4.5pt}{\includegraphics[page=17]{skein}}$, $\raisebox{-4.5pt}{\includegraphics[page=17]{skein}}$, and $\raisebox{-4.5pt}{\includegraphics[page=17]{skein}}$ belong to the same torsion class then the normalized one variable twisted Alexander function satisfies the skein relation
$$\nabla \big(\raisebox{-4.5pt}{\includegraphics[page=16]{skein}}\big) -\nabla\big(\raisebox{-4.5pt}{\includegraphics[page=17]{skein}}\big)=(t-t^{-1})\nabla \big(\raisebox{-4.5pt}{\includegraphics[page=18]{skein}}\big).$$.

\end{theorem}



\section{Examples} \label{sec:examples}

We finish by presenting some explicit calculations of difficult cases of links in $L(p,1)$ where the mentioned invariants fail to detect inequivalent links. The knot notations are taken from the lens space knot table constructed in \cite{gabr}. The Kauffman bracket skein modules and HOMFLY-PT skein modules (evaluated in the standard basis) were computed by the C++ program available in~\cite{code} (the algorithm itself is presented \cite{gabr}).
The Alexander polynomials were computed using SnapPy and SageMath and applying equation~(\ref{eq:alexlpq}). The Kauffman skein modules and Dubrovnik skein modules were computed by hand (for the solid torus and by linearity substituting the solid torus generators with the lens space generators).

\ifx false

\begin{example}

Consider the knots $\ov{4_3}$ and $4_{11}$ in Figure~\ref{fig:ex1}.
We have the following values of the Kauffman bracket skein modules, the HOMFLYPT skein modules and the Alexander polynomial:

\begin{figure}[ht]
	\centering
	\subfigure[$\ov{4_3}$]{\begin{overpic}[page=1]{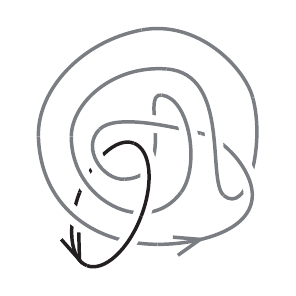}
		\put(15,12){\scriptsize{$p$}}
	\end{overpic}}
	\qquad\qquad\qquad
	\subfigure[$4_{11}$]{\begin{overpic}[page=2]{knots}
		\put(14,11){\scriptsize{$p$}}
	\end{overpic}}
	\caption{Knots $\ov{4_3}$ and $4_{11}$.}\label{fig:ex1}
\end{figure}

\vspace{1em}
$\kbsm(\ov{4_3}) = \begin{cases}
\mej{\left(A^{19}-2 A^{15}+A^{11}-A^7+A^3-A^{-1}\right) x},&
p=2\\
\mej{-A^{12}+A^8+A^4-1+\left(A^{19}-2 A^{15}+A^{11}-A^3\right) x},&
p=3\\
\mej{\left(-A^{16}+2 A^{12}+A^{10}-A^8-2 A^6+A^2+1\right) x},&
p=4\\
\mej{-A^{11}+A^7+A^3-A^{-1}+\left(-A^{16}+2 A^{12}-A^8+1\right)x} &\\
\mej{\quad+\left(A^{11}-2 A^7+A^3\right) x^2}, &
p=5\\
\mej{\left(A^{13}-2 A^5\right) x+\left(-A^9+2 A^5-A\right) x^3}, &
p \geq 6.\\
\end{cases}$

\vspace{1em}
$\kbsm(4_{11}) = \begin{cases}
\mej{\left(A^{19}-2 A^{15}+A^{11}-A^7+A^3-A^{-1}\right) x}
,& p=2\\
\mej{-A^{12}-1+\left(A^{19}-2 A^{15}+A^{11}\right) x}
,& p=3\\
\mej{\left(A^{13}-2 A^{9}-A^7+A^5+A^3-A^{-1}\right) x}
,& p=4\\
\mej{-A^{12}-A^{-1}+\left(-A^{16}+2 A^{12}-A^8\right) x} &\\
\mej{\quad+\left(A^{11}-A^7+A^3\right) x^2}
,& p=5\\
\mej{\left(-A^{16}-A^4-1\right) x+\left(A^{12}-A^8+A^4\right) x^3}
,& p \geq 6.\\
\end{cases}$

\vspace{1em}
$\hsm(\ov{4_3}) = \begin{cases}
\mej{ v^{2}z^{2}t_{1}^{3} + (- z^{4} - 2z^{2} + 2v^{-2}z^{2} + 2v^{-2} - v^{-4})t_{1}}
,&p=2\\
\mej{v^{-2}z^{2}t_{-1}^{2}t_{1} + (2v^{-2} - v^{-4})t_{-1}}
,&p \geq 3.\\
\end{cases}$

\vspace{1em}
$\hsm(4_{11}) = \begin{cases}
\mej{- z^{2}t_{1}^{3} + (z^{2} + 2 - v^{-2}z^{2} - 2v^{-2} + v^{-4})t_{1}}
,& p=2\\
\mej{(- v^{-1}z + v^{-3}z)t_{-1}t_{1} - z^{2}t_{1}^{3} - vz} &\\
\mej{\quad - vz^{-1} + 2v^{-1}z + v^{-1}z^{-1} - v^{-3}z}
,& p=3\\
\mej{t_{-1} - z^{2}t_{1}^{3} +(- vz + v^{-1}z)t_{1}t_{2}}
,& p\geq4.\\
\end{cases}$

\vspace{1em}
$\Delta(\ov{4_3})=
1-2 t+t^2-2 t^p+t^{2 p}+3 t^{p+1}-2 t^{p+2}-2 t^{2 p+1}+t^{2 p+2}.
$

\vspace{1em}
$\Delta(4_{11})=
-1 + 2 t + 2 t^p - 3 t^{p+1} + 2 t^{p+2} + 2 t^{2 p+1} - t^{ 2 p+2}.
$
\vspace{1em}

As the HOMFLYPT skein module and the Alexander polynomial suggest, the links are inequivalent for any $p$, but the Kauffman bracket skein module takes the same values for $p=2$.

\end{example}
\fi

\begin{example}

Consider the knots $5_{76}$ and $\ov{5_{76}}$ in Figure~\ref{fig:ex2}.
The knot $\ov{5_{76}}$ differs from $5_{76}$ by exchanging the crossing on the moving component, which can be interpreted as $\ov{5_{76}}$ being the mirror image of $5_{76}$ under the self-homeomorphism of $T$ that reverses the orientation of the meridian but keeps the orientation of the longitude.
Amphichirality of $5_{76}$ is not detected by the Kauffman bracket skein module for any value of $p$, but detected by the other skein modules and the Alexander polynomial.

\begin{figure}[ht]
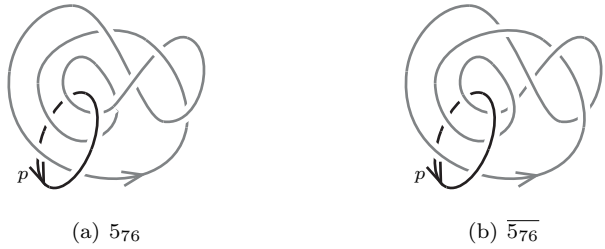

	\centering
	\subfigure[$5_{76}$]{\begin{overpic}[page=3]{knots}
		\put(11,12){\scriptsize{$p$}}
	\end{overpic}}
	\qquad\qquad\qquad
	\subfigure[$\ov{5_{76}}$]{\begin{overpic}[page=9]{knots}
		\put(11,12){\scriptsize{$p$}}
	\end{overpic}}
	\caption{The knot $5_{76}$ and its mirror image.}\label{fig:ex2}
\end{figure}


\vspace{1em}
\begin{tabularx}{0.85\textwidth}{|c|X|}  \hline
$p$ & $\kbsm(5_{76}) = \kbsm(\ov{5_{76}})$
\\  \hline
$2$ & \mej{
x
} \\  \arrayrulecolor{gray}\hline\arrayrulecolor{black}
$3$ & \mej{
A^{13}+A+x(-A^8 + A^4)
} \\  \arrayrulecolor{gray}\hline\arrayrulecolor{black}
$4$ & \mej{
x(-A^{11}+A^7+A^5-A^3-A) 
} \\  \arrayrulecolor{gray}\hline\arrayrulecolor{black}
$\geq\!5$ & \mej{
x(2 A^7-A^3+2A^{-1}) + x^3 (-A^7+A^3-A^{-1}) 
} \\  \hline
\end{tabularx}


\vspace{1em}
\begin{tabularx}{0.85\textwidth}{|c|X|} \hline
$p$ & $\hsm(5_{76})$
\\  \hline
$2$ & \mej{
t_1
} \\  \arrayrulecolor{gray}\hline\arrayrulecolor{black}
$3$ & \mej{
 - v^{-1}z_{-1} + v^{-3}z^{-1} -t_{-1}t_1 \, v^{-1}z
} \\  \arrayrulecolor{gray}\hline\arrayrulecolor{black}
$4$ & \mej{
t_{-1}\,(-z^2 + v^{-2}) - t_1t_2 \,vz 
} \\  \arrayrulecolor{gray}\hline\arrayrulecolor{black}
$5$ & \mej{
-t_{-1}t_{-2}\,vz +t_2\,(- z^2 + v^{-2})
} \\  \arrayrulecolor{gray}\hline\arrayrulecolor{black}
$\geq\!6$ & \mej{
-t_{1}t_{2}\,  vz+t_{3}\,(- z^{2} + v^{-2})
} \\  \hline
\end{tabularx}

\vspace{1em}
\begin{tabularx}{0.85\textwidth}{|c|X|}  \hline
$p$ & $\hsm(\ov{5_{76}})$
\\  \hline
$2$ & \mej{
t_{1} \,(2v^{-2}z^{2} + 3v^{-2} - 2v^{-4}z^{4} - 4v^{-4}z^{2} - 3v^{-4} + 2v^{-6}z^{2} + v^{-6}) 
+ t_{1}^{3}\,(z^{2} - v^{-2}z^{4} - v^{-2}z^{2})
} \\  \arrayrulecolor{gray}\hline\arrayrulecolor{black}
$3$ & \mej{
- 2v^{-1}z - v^{-1}z^{-1} + 2v^{-3}z^{3} + 4v^{-3}z + v^{-3}z^{-1} - 2v^{-5}z^{3} - 2v^{-5}z }\\&\mej{
+t_{-1}^{3} \,(z^{2} - v^{-2}z^{4} - v^{-2}z^{2})
+t_{-1}t_{1}\,(- 2v^{-3}z + 2v^{-5}z^{3} + v^{-5}z)
} \\  \arrayrulecolor{gray}\hline\arrayrulecolor{black}
$4$ & \mej{
t_{-1}^{3} \, (z^{2} - v^{-2}z^{4} - v^{-2}z^{2})
+ t_{1}\,  (2v^{-2}z^{2} + v^{-2} - 2v^{-4}z^{4} - 2v^{-4}z^{2}) }\\&\mej{
+ t_{-1}t_{2}\,(- 2v^{-3}z+ 2v^{-5}z^{3} + v^{-5}z) 
} \\  \arrayrulecolor{gray}\hline\arrayrulecolor{black}
$5$ & \mej{
t_{-2}t_{-1}\, (- 2v^{-1}z + 2v^{-3}z^{3} + v^{-3}z)
+t_{-1}^{3}\,( z^{2} - v^{-2}z^{4}- v^{-2}z^{2}) 
+t_{2}\, (v^{-2} - v^{-4}z^{2})
} \\  \arrayrulecolor{gray}\hline\arrayrulecolor{black}
$\geq\!6$ & \mej{
t_{-3}\,(v^{-2} - v^{-4}z^{2}) 
+t_{-2}t_{-1}\,(- 2v^{-1}z + 2v^{-3}z^{3} + v^{-3}z)
 + t_{-1}^{3}\,(z^{2} - v^{-2}z^{4} - v^{-2}z^{2})
} \\   \hline
\end{tabularx}


\vspace{1em}
\begin{tabularx}{0.85\textwidth}{|c|X|} \hline
$p$ & $\ksm(5_{76})$
\\ \hline
$2$ & \mej{
t_1\, (-z+a^{-1}+z a^2)
} \\  \arrayrulecolor{gray}\hline\arrayrulecolor{black}
$3$ & \mej{
z a^2+a^3+z a^4-a^2 z^{-1}-a^4 z^{-1}-a^3 z^2 
+ t_1 \,(-z-a z^2+a^{-1} z^2+z^3-a^2 z^3)  }\\&\mej{
+ t_1^2\,(z a^2+a^3 z^2)
} \\  \arrayrulecolor{gray}\hline\arrayrulecolor{black}
$4$ & \mej{
t_1\, (-z-a^5+a^{-1} z^2+a^3 z^2+a^5 z^2+z^3+a^4 z^3) 
+ t_1 t_2 \,(-z-a z^2)
} \\  \arrayrulecolor{gray}\hline\arrayrulecolor{black}
$5$ & \mej{
t_1 \,(-z+a^{-1} z^2+z^3)
+ t_2 \,(-a^7+a^5 z^2+a^7 z^2+a^6 z^3)
+ t_1 t_2 \,(-z-a z^2)
} \\  \arrayrulecolor{gray}\hline\arrayrulecolor{black}
$\geq\!6$ & \mej{
t_1\, (-z+a^{-1} z^2+z^3)
+ t_1 t_2 \,(-z-a z^2) 
+ t_3\, (-a^3+a z^2+a^3 z^2+a^2 z^3)
} \\  \hline
\end{tabularx}


\vspace{1em}
\begin{tabularx}{0.85\textwidth}{|c|X|}\hline
$p$ & $\ksm(\ov{5_{76}})$
\\  \hline
$2$ & \mej{
t_1\, (-3+2 a z-2 z a^{-1}-2 a^2-2 z a^2-a^3-a^5+8 z^2+5 a^2 z^2+4 a^3 z^2+2 a^4 z^2 }\\&\mej{
+3 a z^3+2 a^{-1} z^3+3 a^2 z^3-3 z^4-7 a^2 z^4-5 a z^5) 
+ t_1^2\, (-a^3 z^2-a^2 z^3)  }\\&\mej{
+ t_1^3\, (-z^2-a^3 z^2-a^2 z^3+z^4+a^2 z^4+a z^5)
} \\  \arrayrulecolor{gray}\hline\arrayrulecolor{black}
$3$ & \mej{
3 a z+a^2+4 z a^3+z a^4+z a^5+z a^6-a z^{-1}-a^3 z^{-1}-2 a^2 z^2+a^3 z^2-2 a z^3-5 a^3 z^3 }\\&\mej{
-a^4 z^3-2 a^5 z^3 +2 a^3 z^5
+ t_1 \,(-2 z a^{-1}-z a^2+z^2+2 a^2 z^2+a^3 z^2 +a^5 z^2+2 a z^3}\\&\mej{
+2 a^{-1} z^3+a^2 z^3+a^4 z^3+z^4-3 a^2 z^4-2 a^4 z^4-a z^5-2 a^3 z^5) }\\&\mej{
+ t_1^2 \,(-2 z a^3-z a^6-a^3 z^2-a^5 z^2-a^2 z^3+2 a^3 z^3+2 a^5 z^3+2 a^4 z^4) }\\&\mej{
+ t_1^3\, (-z^2-a^3 z^2-a^2 z^3+z^4+a^2 z^4+a z^5)
} \\  \arrayrulecolor{gray}\hline\arrayrulecolor{black}
$4$ & \mej{
 t_1 \,(-2 z a^{-1}-z a^2-a^4+z^2+a^3 z^2+a^4 z^2+a^6 z^2+2 a z^3+2 a^{-1} z^3+a^2 z^3+a^5 z^3 }\\&\mej{
  +z^4-a^2 z^4-a z^5)
+ t_1^2\, (-a^3 z^2-a^2 z^3) 
+ t_1 t_2\, (2 a z+z a^4+a^3 z^2-2 a z^3-2 a^3 z^3 }\\&\mej{
-2 a^2 z^4)
+ t_1^3\, (-z^2-a^3 z^2-a^2 z^3+z^4+a^2 z^4+a z^5)
} \\  \arrayrulecolor{gray}\hline\arrayrulecolor{black}
$5$ & \mej{
t_1\, (-2 z a^{-1}-z a^2+z^2+a^3 z^2+2 a z^3+2 a^{-1} z^3+a^2 z^3+z^4-a^2 z^4-a z^5) }\\&\mej{
+ t_1^2 \,(-a^3 z^2-a^2 z^3)
+ t_2 \,(-a^6+a^6 z^2+a^8 z^2+a^7 z^3)
+ t_1 t_2\, (2 a z+z a^4+a^3 z^2  }\\&\mej{
-2 a z^3-2 a^3 z^3-2 a^2 z^4)
+ t_1^3\, (-z^2-a^3 z^2-a^2 z^3+z^4+a^2 z^4+a z^5)
} \\  \arrayrulecolor{gray}\hline\arrayrulecolor{black}
$\geq\!6$ & \mej{
 t_1 \,(-2 z a^{-1}-z a^2+z^2+a^3 z^2+2 a z^3+2 a^{-1} z^3+a^2 z^3+z^4-a^2 z^4-a z^5) }\\&\mej{
+ t_1^2 \,(-a^3 z^2-a^2 z^3) 
+ t_1^3 \,(-z^2-a^3 z^2-a^2 z^3+z^4+a^2 z^4+a z^5)}\\&\mej{
+ t_1 t_2 \,(2 a z+z a^4+a^3 z^2 
-2 a z^3-2 a^3 z^3-2 a^2 z^4)
+ t_3 \,(-a^2+a^2 z^2+a^4 z^2+a^3 z^3)
} \\  \hline
\end{tabularx}


\vspace{1em}
\begin{tabularx}{0.85\textwidth}{|c|X|}\hline
$p$ & $\dsm(5_{76})$
\\  \hline
$2$ & \mej{
t_1\, (z+a^{-1}-z a^2)
} \\  \arrayrulecolor{gray}\hline\arrayrulecolor{black}
$3$ & \mej{
-z a^2+a^3+z a^4-a^2 z^{-1}+a^4 z^{-1}+a^3 z^2
+ t_1\, (z-a z^2+a^{-1} z^2+z^3-a^2 z^3) }\\&\mej{
+ t_1^2 \,(-z a^2-a^3 z^2)
} \\  \arrayrulecolor{gray}\hline\arrayrulecolor{black}
$4$ & \mej{
 t_1 \,(z+a^5+a^{-1} z^2-a^3 z^2+a^5 z^2+z^3-a^4 z^3)
+ t_1 t_2 \,(-z-a z^2)
} \\  \arrayrulecolor{gray}\hline\arrayrulecolor{black}
$5$ & \mej{
 t_1 \,(z+a^{-1} z^2+z^3)
+ t_2 \,(a^7-a^5 z^2+a^7 z^2-a^6 z^3)
+ t_1 t_2\, (-z-a z^2)
} \\  \arrayrulecolor{gray}\hline\arrayrulecolor{black}
$\geq\!6$ & \mej{
 t_1 \,(z+a^{-1} z^2+z^3)
+ t_1 t_2 \,(-z-a z^2)
+ t_3\, (a^3-a z^2+a^3 z^2-a^2 z^3)
} \\ \hline
\end{tabularx}

\vspace{1em}
\begin{tabularx}{0.85\textwidth}{|c|X|}  \hline
$p$ & $\dsm(\ov{5_{76}})$
\\ \hline
$2$ & \mej{
t_1 \,(3-2 a z+2 z a^{-1}-2 a^2-a^3+a^5+2 z^2-5 a^2 z^2+2 a^4 z^2-a z^3+2 a^{-1} z^3-a^2 z^3 }\\&\mej{
+z^4-a^2 z^4-a z^5)
+t_1^2\, (a^3 z^2-a^2 z^3)
+ t_1^3\, (z^2-a^3 z^2+a^2 z^3+z^4-a^2 z^4+a z^5)
} \\  \arrayrulecolor{gray}\hline\arrayrulecolor{black}
$3$ & \mej{
-3 a z+a^2+4 z a^3+z a^4-z a^5-z a^6-a z^{-1}+a^3 z^{-1}+2 a^2 z^2-a^3 z^2-2 a z^3 }\\&\mej{
+5 a^3 z^3+a^4 z^3-2 a^5 z^3
+ t_1\, (2 z a^{-1}-z a^2-z^2-2 a^2 z^2+a^3 z^2+a^5 z^2-2 a z^3}\\&\mej{
+2 a^{-1} z^3  -a^2 z^3-a^4 z^3+z^4-a^2 z^4+2 a^4 z^4-a z^5-2 a^3 z^5) }\\&\mej{
+ t_1^2\, (-2 z a^3+z a^6+a^3 z^2-a^5 z^2-a^2 z^3-2 a^3 z^3+2 a^5 z^3-2 a^4 z^4)+2 a^3 z^5 }\\&\mej{
+ t_1^3 \,(z^2-a^3 z^2+a^2 z^3+z^4-a^2 z^4+a z^5)
} \\  \arrayrulecolor{gray}\hline\arrayrulecolor{black}
$4$ & \mej{
t_1\, (2 z a^{-1}-z a^2+a^4-z^2+a^3 z^2+a^4 z^2-a^6 z^2-2 a z^3+2 a^{-1} z^3-a^2 z^3+a^5 z^3 }\\&\mej{
+z^4+a^2 z^4-a z^5)
+ t_1^2 \,(a^3 z^2-a^2 z^3)
+ t_1 t_2 \,(-2 a z+z a^4-a^3 z^2-2 a z^3+2 a^3 z^3 }\\&\mej{
 -2 a^2 z^4)
+ t_1^3 \,(z^2-a^3 z^2+a^2 z^3+z^4-a^2 z^4+a z^5)
} \\  \arrayrulecolor{gray}\hline\arrayrulecolor{black}
$5$ & \mej{
t_1\, (2 z a^{-1}-z a^2-z^2+a^3 z^2-2 a z^3+2 a^{-1} z^3-a^2 z^3+z^4+a^2 z^4-a z^5) }\\&\mej{
+ t_1^2 \,(a^3 z^2-a^2 z^3)
+ t_2\, (a^6+a^6 z^2-a^8 z^2+a^7 z^3) 
+ t_1 t_2 \,(-2 a z+z a^4 -a^3 z^2-2 a z^3}\\&\mej{
+2 a^3 z^3-2 a^2 z^4)
+ t_1^3\, (z^2-a^3 z^2+a^2 z^3+z^4-a^2 z^4+a z^5)
} \\  \arrayrulecolor{gray}\hline\arrayrulecolor{black}
$\geq\!6$ & \mej{
t_1 \,(2 z a^{-1}-z a^2-z^2+a^3 z^2-2 a z^3+2 a^{-1} z^3-a^2 z^3+z^4+a^2 z^4-a z^5) }\\&\mej{
+ t_1^2\, (a^3 z^2-a^2 z^3)
+ t_1^3\, (z^2-a^3 z^2+a^2 z^3+z^4-a^2 z^4+a z^5) }\\&\mej{
+ t_1 t_2 \,(-2 a z+z a^4-a^3 z^2-2 a z^3+2 a^3 z^3-2 a^2 z^4)
+ t_3 \,(a^2+a^2 z^2-a^4 z^2+a^3 z^3)
} \\   \hline
\end{tabularx}

\vspace{1em}
$\Delta(5_{76})  = 
-t^{2 p-1}-2 t^{3 p-2}+t^{4 p-2}-2 t^p+1.$

\vspace{1em}
$\Delta(\ov{5_{76}}) = 
-t^{2 p+1}-2 t^{3 p+2}+t^{4 p+2}-2 t^p+1.
$

\end{example}

\begin{example}

The knots $5_{26}$ and $5_{27}$ in Figure~\ref{fig:ex3} differ by exchanging both the orientation of the fixed and mixed sublinks, which can be interpreted as $5_{27}$ being the image of $5_{26}$ under the self-homomorphism of the torus $T$ that reverses both the meridian and the longitude (a so-called \emph{flip} in the language of~\cite{BM}, see also~\cite{cm1}). The question whether $5_{26} \neq 5_{27}$ is equivalent to the question whether the links are non-invertible.

\begin{figure}[ht]
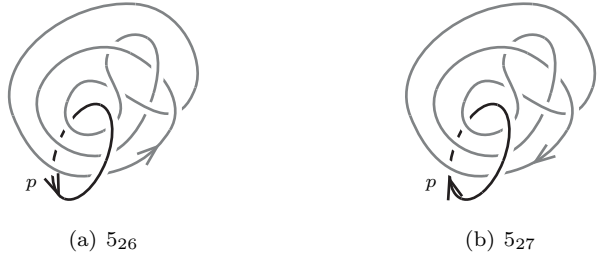

	\centering
	\subfigure[$5_{26}$]{\begin{overpic}[page=5]{knots}
		\put(16,12){\scriptsize{$p$}}
	\end{overpic}}
	\qquad\qquad\qquad
	\subfigure[$5_{27}$]{\begin{overpic}[page=10]{knots}
		\put(17,12){\scriptsize{$p$}}
	\end{overpic}}
	\caption{The knots $5_{26}$ and its flip $5_{27}$.}\label{fig:ex3}
\end{figure}

Non-invertible links were studied by Whitten~\cite{whitten} and are hard to detect, although in the case when the links are hyperbolic (most are), modern computational techniques using canonical triangulations of the link complements enable us to verifiably recognize them~\cite{W}. 

It is shown in~\cite{gabr} that $5_{26}$ and $5_{27}$ are non-isotopic in any lens space $L(p,1)$, but due to the symmetric nature of the two knots, none of our invariants are able to detect this.


\vspace{1em}
\begin{tabularx}{0.85\textwidth}{|c|X|} \hline
$p$ & $\kbsm(5_{26}) = \kbsm(5_{27})$
\\  \hline
$2$ & \mej{
x\,(-A^{24}+3 A^{20}-2 A^{16}+3 A^{12}
-3 A^8+2 A^4-1) 
} \\  \arrayrulecolor{gray}\hline\arrayrulecolor{black}
$3$ & \mej{
A^{17}-A^{13}+A^5-A
- x \,A^8 (A^{16} -3 A^{12} +2 A^8-2 A^4+1)
} \\  \arrayrulecolor{gray}\hline\arrayrulecolor{black}
$4$ & \mej{
-x\,(A^{18}-3 A^{14}-A^{12} + 2 A^{10}+2 A^8-2 A^6-2 A^4+A^2+1)
} \\  \arrayrulecolor{gray}\hline\arrayrulecolor{black}
$5$ & \mej{
A^{16}-A^{12}+A^4-1
+x \,(A^{21}-3 A^{17}+2A^{13}-2 A^9+A^5) }\\&\mej{
+x^2 \,(-A^{16}+2 A^{12}-2 A^8+A^4)
} \\  \arrayrulecolor{gray}\hline\arrayrulecolor{black}
$\geq\!6$ & \mej{
x \,(A^{21}-A^{17}-A^{13}+A^5-A) 
+ x^3\,(-A^{17}+2 A^{13}-2 A^9+A^5) 
} \\   \hline
\end{tabularx}


\vspace{1em}
\begin{tabularx}{0.85\textwidth}{|c|X|} \hline
$p$ & $\hsm(5_{26}) = \hsm(5_{27})$
\\  \hline
$2$ & \mej{
t_{1}^{3} \, (- v^{6}z^{4}- v^{6}z^{2})
+ t_{1}\,(v^{4}z^{6}+ 3v^{4}z^{4} + 2v^{4}z^{2} - v^{4} - v^{2}z^{4} + 2v^{2})
} \\  \arrayrulecolor{gray}\hline\arrayrulecolor{black}
$3$ & \mej{
t_{-1}^{2}\,(v^{3}z^{3} + v^{3}z) 
+t_{-1}t_{1}^{2}\,(- v^{2}z^{4} - v^{2}z^{2})
+ t_{1}\,(v^{4}z^{2} + v^{2}z^{4} + v^{2}z^{2} + v^{2})
} \\  \arrayrulecolor{gray}\hline\arrayrulecolor{black}
$\geq\!4$ & \mej{
t_{-1}t_{1}^{2}\, (- v^{2}z^{4} - v^{2}z^{2})
+ t_{-1}t_{2}\,(vz^{3} + vz)
+ t_{1}\, (v^{4}z^{2} + v^{2})
} \\ \hline
\end{tabularx}


\vspace{1em}
\begin{tabularx}{0.85\textwidth}{|c|X|} \hline
$p$ & $\ksm(5_{26}) =\ksm(5_{27})$
\\  \hline
$2$ & \mej{
-a z^2-a^{-1} z^2 +2 z^3+a^2 z^3-a z^4
+ t_1\, (-2 a-z+2 a z-z a^{-2}+z a^2-a^3+z a^3 }\\&\mej{
-2 a z^2-2 a^{-1} z^2-2 a^2 z^2+z^3-2 a z^3+2 a^{-2} z^3-a^2 z^3+2 a z^4+4 a^{-1} z^4+2 z^5) }\\&\mej{
+ t_1^2\, (-a z^2+a^2 z^3+a z^4)
-t_1^3\, a  z^2
} \\  \arrayrulecolor{gray}\hline\arrayrulecolor{black}
$3$ & \mej{
-z+a^2-z a^2-z a^3+a^4-a^2 z^2-a^3 z^2-a^4 z^2 +2 z^3+3 a^2 z^3+a^3 z^3-a z^4+a^3 z^4 }\\&\mej{
-a^2 z^5
+ t_1\, (-a-a z^2-a^2 z^2+a^3 z^2+2 z^3-a^3 z^3+3 a z^4+a^3 z^4+a^2 z^5) }\\&\mej{
+ t_1^2\, (-z a^4-2 a z^2+a^4 z^2-a^2 z^3+a z^4-a^3 z^4)
-t_1^3\, a z^2
} \\  \arrayrulecolor{gray}\hline\arrayrulecolor{black}
$4$ & \mej{
-a z^2-a^{-1} z^2+z^3
+ t_1 \,(-a+z a^5-a z^2-a^2 z^2-a^3 z^2-a^5 z^2+z^3-a^4 z^3+a z^4) }\\&\mej{
+ t_1^2 \,(-a z^2+a^2 z^3+a z^4)
+ t_2\, a z^4
- t_1^3\, a z^2
+ t_1 t_2\, (z a^2+a^{-1} z^2-a^2 z^2+2 z^3+a z^4)
} \\  \arrayrulecolor{gray}\hline\arrayrulecolor{black}
$5$ & \mej{
-a z^2-a^{-1} z^2 +z^3
+ t_1 \,(-a-a z^2-a^2 z^2+z^3+a z^4)
+ t_1^2 \,a(- z^2+a z^3+ z^4) 
+ t_2 }\\&\mej{ \,(z a^7-a^5 z^2-a^7 z^2-a^6 z^3+a z^4)
-t_1^3\, a  z^2 
+ t_1 t_2 (z a^2+a^{-1} z^2-a^2 z^2+2 z^3+a z^4)
} \\  \arrayrulecolor{gray}\hline\arrayrulecolor{black}
$\geq\!6$ & \mej{
-a z^2- z^2/a +z^3
+ t_1\, (-a-a z^2-a^2 z^2+z^3+a z^4)
+ t_1^2 a(- z^2+a z^3+ z^4)
+ t_2 a  z^4 }\\&\mej{
-t_1^3\, a z^2
+ t_1 t_2\, (z a^2+a^{-1} z^2-a^2 z^2+2 z^3+a z^4)
+ t_3\, (z a^3-a z^2-a^3 z^2-a^2 z^3)
} \\   \hline
\end{tabularx}


\vspace{1em}
\begin{tabularx}{0.85\textwidth}{|c|X|}  \hline
$p$ &  $\dsm(5_{26}) =\dsm(5_{27})$
\\  \hline
$2$ & \mej{
-a z^2+a^{-1} z^2 -2 z^3+a^2 z^3+a z^4
+ t_1 \,(z+2 a z-z a^{-2}+z a^2+a^3-z a^3+2 a z^2 }\\&\mej{
-2 a^{-1} z^2-z^3+a^2 z^3)
+ t_1^2 \,(a z^2+a^2 z^3-a z^4)
- t_1^3 \,a  z^2
} \\  \arrayrulecolor{gray}\hline\arrayrulecolor{black}
$3$ & \mej{
z-a^2-z a^2+z a^3+a^4-2 a z^2+2 a^{-1} z^2-a^2 z^2-a^3 z^2+a^4 z^2-3 a^2 z^3+a^3 z^3 }\\&\mej{
-a z^4-a^3 z^4-a^2 z^5
+ t_1 \,(a-a z^2+a^2 z^2+a^3 z^2-a^3 z^3+a z^4+a^3 z^4+a^2 z^5) }\\&\mej{
+ t_1^2 \,(z a^4+2 a z^2-a^4 z^2+3 a^2 z^3-a z^4+a^3 z^4)
-t_1^3 \,a  z^2
} \\  \arrayrulecolor{gray}\hline\arrayrulecolor{black}
$4$ & \mej{
-a z^2+a^{-1} z^2-z^3
+ t_1 \,(a+z a^5-a z^2+a^2 z^2-a^3 z^2+a^5 z^2-z^3-a^4 z^3-a z^4) }\\&\mej{
+ t_1^2 \,(a z^2+a^2 z^3-a z^4)
+t_2 \,a z^4
+ t_1 t_2 \,(z a^2+a^{-1} z^2-a^2 z^2+2 z^3+a z^4)
-t_1^3\, a  z^2
} \\  \arrayrulecolor{gray}\hline\arrayrulecolor{black}
$5$ & \mej{
-a z^2+a^{-1} z^2 -z^3
-t_1^3 \,a z^2
+ t_1\, (a-a z^2+a^2 z^2-z^3-a z^4)
+ t_1^2 \,a(z^2+a z^3- z^4) }\\&\mej{
+ t_1 t_2\, (z a^2+a^{-1} z^2-a^2 z^2+2 z^3+a z^4)
+ t_2 \,(z a^7-a^5 z^2+a^7 z^2-a^6 z^3+a z^4)
} \\  \arrayrulecolor{gray}\hline\arrayrulecolor{black}
$\geq\!6$ & \mej{
-a z^2+a^{-1} z^2 -z^3
+ t_1\, (a-a z^2+a^2 z^2-z^3-a z^4)
+ t_1^2 \,a( z^2+a z^3-z^4)
+ t_2 \,a z^4 }\\&\mej{
-t_1^3\, a z^2
+ t_1 t_2\, (z a^2+a^{-1} z^2-a^2 z^2+2 z^3+a z^4)
+ t_3\, (z a^3-a z^2+a^3 z^2-a^2 z^3)
} \\   \hline
\end{tabularx}

\vspace{1em}
$\Delta(5_{26}) =\Delta(5_{27}) = (t^2-t+1) (t^{p+1}-t+1) (t^{p+1}-t^p+1).$

\end{example}

\subsubsection*{Acknowledgments}
The first author was supported by the Slovenian Research Agency grants J1-8131, J1-7025, and N1-0064.


\end{document}